\documentclass[11pt]{article}

\usepackage[margin=1in]{geometry}
\usepackage{setspace}
\usepackage{amsmath, amssymb, amsthm, mathtools}
\usepackage{bm}
\usepackage{authblk}
\usepackage{orcidlink}

\usepackage{natbib}
\usepackage{hyperref}
\hypersetup{
	colorlinks=true,
	linkcolor=blue,
	citecolor=blue,
	urlcolor=blue
}

\theoremstyle{plain}
\newtheorem{theorem}{Theorem}
\newtheorem{proposition}{Proposition}
\newtheorem{lemma}{Lemma}
\newtheorem{corollary}{Corollary}

\theoremstyle{definition}
\newtheorem{definition}{Definition}

\theoremstyle{remark}
\newtheorem{remark}{Remark}

\newcommand{\E}{\mathbb{E}}

\newcommand{\KL}{\mathrm{KL}}

\newcommand{\diag}{\mathrm{diag}}

\title{Bayesian Effective Dimension from Information Growth}

\author{Sayantan Banerjee \orcidlink{0000-0001-5414-4817}}
\affil{Operations Management \& Quantitative Techniques Area\\Indian Institute of Management Indore \\ Rau-Pithampur Road, Indore MP 453 556
}

\date{Email: \texttt{sayantanb@iimidr.ac.in}}

\begin{document}
	\maketitle
	
\begin{abstract}
	Nominal dimension can greatly overstate how much a Bayesian model actually learns from data. We quantify this gap through a prior-dependent effective dimension based on parameter--data mutual information, normalized so that regular parametric information growth provides the benchmark scale. The main theoretical result shows that, in infinite-dimensional Gaussian experiments, information growth is determined by the spectral counting function: regularly varying spectra yield a sharp asymptotic law with an explicit constant, while exponential spectral decay produces a distinct polylogarithmic regime. In Gaussian regression the same information scale has a spectral representation and is linked exactly to ridge degrees of freedom. For scalar Gaussian scale-mixture priors, we separate information about latent scales from conditional Gaussian information and show that the full marginal information remains finite under polynomially heavy-tailed mixing.
\end{abstract}

	
	
	\sloppy

\section{Introduction}
\label{sec:introduction}

Dimension has a less straightforward meaning in Bayesian problems than the number of coordinates in a parameter vector might suggest. A regression model may contain thousands of coefficients while the data provide substantial information about only a smaller collection of directions. In an inverse problem, the parameter may live in an infinite-dimensional space even though only finitely many modes are appreciably resolved at a given noise level. Shrinkage and regularization create related effects by allocating learning unevenly across coordinates or eigendirections. These phenomena are familiar, although they appear under different names: sparsity, degrees of freedom, intrinsic dimension, effective rank, or regularization.

Several established theories formalize different aspects of this idea. Degrees of freedom quantify the flexibility of fitted procedures and have familiar spectral forms in regularized regression \citep{efron2004least,hastie2009elements}. Bayesian effective-parameter constructions arise in evidence-based regularization and posterior or predictive complexity measures \citep{mackay1992bayesian,spiegelhalter2002bayesian}. Closely related spectral quantities appear in kernel regularization, where effective dimension is expressed through traces of regularized covariance operators and enters learning-rate theory \citep{caponnetto2007optimal}. In inverse problems, intrinsic-dimension constructions describe the interaction between an infinite-dimensional prior and the likelihood, often with computational consequences \citep{agapiou2017importance}. Singular models provide another notion of effective dimension through learning coefficients or real log canonical thresholds, where the regular coefficient \(d/2\) in Bayesian asymptotics can be replaced by a smaller geometric quantity \citep{watanabe2009algebraic}. Sparsity offers another route to reduced effective complexity by restricting attention to a comparatively small active subset of coordinates, although its organizing principle is support structure rather than information accumulation. These objects answer related but distinct questions.

Parameter--data dependence provides an older Bayesian route to complexity. If \(\Theta\sim\Pi\) and \(X^{(n)}\) is generated conditionally on \(\Theta\), then
\[
I(\Theta;X^{(n)})
=
\E\!\left[
\KL\{\Pi(\cdot\mid X^{(n)})\|\Pi\}
\right],
\]
so mutual information is the expected prior-to-posterior Kullback--Leibler gain. This interpretation goes back at least to Lindley's measure of the information supplied by an experiment \citep{lindley1956measure}, and subsequently became central in Bayesian design and decision theory \citep{bernardo1979expected}. Related developments appear in universal coding and redundancy \citep{rissanen1986,clarke1990information}, cumulative Bayesian prediction and relative-entropy risk \citep{haussler1997mutual}, and analyses of the information supplied by a sample about an unknown parameter \citep{ebrahimi2010sample}.

The closest Bayesian precedent is the dependence-based view developed by \citet{van2004association,van2012bayesian}. In that framework Bayesian model complexity is treated as a feature of the joint law of random parameters and observations, and mutual-information-type quantities measure the association generated by Bayesian learning. Gaussian models already lead to log-determinant information expressions, while posterior-predictive dependence can recover familiar effective-parameter and degrees-of-freedom quantities. Thus neither parameter--data mutual information nor its interpretation as Bayesian complexity is new. \citet{handley2019bayesian} pursue a different information-based notion of dimensionality, using fluctuations of posterior Shannon information to quantify the number of constrained directions.

There is also an important asymptotic precedent for interpreting information growth through dimension. In regular \(d\)-dimensional experiments, classical redundancy and Bayesian information calculations give \(I(\Theta;X^{(n)})=(d/2)\log n+O(1)\) under suitable regularity conditions \citep{clarke1990information}. Closely related ideas arise in predictive-information theory, where logarithmic information growth has been interpreted through the number of learnable parameters, with nonparametric models exhibiting different growth regimes \citep{bialek2001predictability}. These results suggest that the \emph{rate} at which information accumulates can itself carry dimensional information.

In this work, we study the specific finite-\(n\) normalization
\[
d_{\mathrm{eff}}(n)
=
\frac{2I(\Theta;X^{(n)})}{\log n},
\]
which expresses parameter--data information in units of the regular-parametric contribution \((1/2)\log n\). For a regular fixed-prior \(d\)-dimensional model, \(d_{\mathrm{eff}}(n)\to d\). In less regular experiments it need not be integer-valued or coincide with the algebraic dimension of the parameter space. The quantity is attached to the full prior--sampling experiment and measures the amount of expected Bayesian learning relative to the regular parametric benchmark.

The finite-\(n\) normalization is useful when information is distributed unevenly across directions. In Gaussian regression, mutual information is a log-determinant determined by the singular values of the design and the prior scale. Weak modes contribute little, while strongly identified modes contribute logarithmically. Infinite-dimensional Gaussian inverse problems have a closely related log-determinant structure: expected posterior-to-prior information and Bayesian \(D\)-optimality have been developed rigorously in Hilbert-space settings \citep{alexanderian2016bayesian}. Our use of this structure is different: the log-determinant is normalized by regular-parametric information growth and interpreted as a dimension scale for a specified prior--sampling experiment.

This spectral viewpoint also yields a precise link with classical ridge degrees of freedom. Dependence-based connections with degrees of freedom have previously been developed in Bayesian model complexity \citep{van2012bayesian}, while kernel regularization uses closely related trace-resolvent quantities as effective dimensions \citep{caponnetto2007optimal}. More generally, the relation between Gaussian mutual information and estimation error is captured by the I--MMSE identity \citep{guo2005mutual}. In the conjugate Gaussian regression setting considered here, parameter--data mutual information satisfies a particularly simple identity: its derivative with respect to log prior-to-noise ratio is exactly the ridge degrees of freedom. Equivalently, mutual information is the accumulated ridge degrees-of-freedom curve over signal-to-noise scale.

In infinite-dimensional Gaussian experiments, the spectral viewpoint leads to a more general information-growth principle. Writing \(N(t)\) for the number of effective signal-to-noise eigenvalues exceeding \(t\), we show that regular variation of \(N(t)\) determines the asymptotic growth of mutual information and hence of the proposed dimension scale. The regular-variation and Tauberian tools used in this argument are classical \citep{bingham1987regular}; the contribution lies in applying them to the parameter--data information functional and its regular-parametric normalization. Polynomial spectral decay yields polynomial information growth with a logarithmic normalization, while exponential decay produces a qualitatively slower polylogarithmic regime.

Spectral information-growth bounds also occur in Gaussian-process learning. In particular, maximal information gain in Gaussian-process bandits can be bounded through kernel eigendecay, with polynomial and faster spectral regimes leading to different information-growth rates \citep{vakili2021information}. The object studied there is maximal information gain over sequentially chosen observations and is used to control regret. Here the prior--sampling experiment is fixed, the mutual information concerns the Bayesian parameter and its observations, and the quantity of interest is its normalization as a dimension scale. The resulting sequence-model results cover canonical mildly and severely ill-posed inverse problems and make explicit how prior regularity and operator ill-posedness enter the information scale.

Hierarchical shrinkage provides a second extension. Global--local priors such as the horseshoe introduce latent scales that modify the conditional signal-to-noise ratios \citep{polson2010shrink,carvalho2010horseshoe,bhadra2019lasso}. Their theory has mainly emphasized estimation, contraction, and adaptation \citep{datta2013asymptotic,van2014horseshoe,van2017adaptive,banerjee2026bayesian}. From the information viewpoint, the hierarchy gives an exact decomposition of marginal parameter--data information into conditional Gaussian information and information learned about the latent scales. Polynomially heavy-tailed scales can have infinite second moments while retaining finite conditional Gaussian information; we further show that the full marginal parameter–data information remains finite under the same tail condition.

The prior dependence of \(d_{\mathrm{eff}}(n)\) is essential. Different priors applied to the same likelihood need not leave the same amount of information for the data to supply, and an \(n\)-dependent prior can alter the normalized growth rate substantially. The information index also forms part of the calibration: sample size is canonical in ordinary repeated-sampling experiments, while inverse noise or precision plays that role in sequence models. Multiplicative re-scaling of the information index does not affect the leading asymptotic dimension, whereas more substantial re-indexing can. We therefore regard \(d_{\mathrm{eff}}(n)\) as a property of a specified sequence of Bayesian experiments.

\subsection*{Contributions}

The paper develops a finite-sample Bayesian effective dimension by calibrating parameter--data mutual information against the regular-parametric benchmark \((1/2)\log n\). We record its basic structural properties and derive exact Gaussian representations, including the identity linking the derivative of parameter--data information along the log signal-to-noise scale to ridge degrees of freedom.

The main theoretical results concern infinite-dimensional Gaussian experiments. We show that the spectral counting function determines information growth, yielding a sharp regularly varying asymptotic with an explicit constant and a distinct polylogarithmic regime under exponential spectral decay. For scalar Gaussian scale-mixture and global--local priors, we decompose parameter--data information into latent-scale and conditional Gaussian components, and show that polynomially heavy-tailed scales may have infinite second moments while the full marginal parameter--data information remains finite.

\subsection*{Organization}

Section~\ref{sec:effective-dimension} defines the effective dimension and develops its regular-model calibration and structural properties. Section~\ref{sec:gaussian-examples} gives the Gaussian spectral calculations and the connection with ridge degrees of freedom. Section~\ref{sec:consequences} develops the spectral counting theory for infinite-dimensional Gaussian experiments and derives polynomial and exponential information-growth regimes. Section~\ref{sec:shrinkage} develops the scale-mixture and global--local calculations. Section~\ref{sec:open-problems} concludes with several open directions.

\subsection*{Notation}

For two positive sequences \(a_n\) and \(b_n\), we write \(a_n\lesssim b_n\) if \(a_n\leq Cb_n\) for some constant \(C>0\), \(a_n\asymp b_n\) if both \(a_n\lesssim b_n\) and \(b_n\lesssim a_n\), and \(a_n\sim b_n\) if \(a_n/b_n\to1\). For sets $A$ and $B$, we write \(A\subset B\) for ordinary set inclusion and \(A\Subset B\) when the closure of \(A\) is compact and contained in the interior of \(B\). The notation \(x_+=\max(x,0)\) is used for the positive part. For probability measures \(P\) and \(Q\), \(\KL(P\|Q)\) denotes Kullback--Leibler divergence. \(I(U;V)\) denotes mutual information. Expectations are written as \(\E\), with subscripts added when the underlying law needs to be emphasized. For a symmetric matrix \(A\), \(\lambda_j(A)\) denotes its eigenvalues and \(\operatorname{tr}(A)\) and \(\det(A)\) its trace and determinant respectively. Throughout, \(n\) indexes the information level of the experiment, typically sample size or inverse noise level.	

\section{Bayesian experiments and effective dimension}
\label{sec:effective-dimension}

\subsection{Bayesian experiments and information gain}
\label{subsec:bayes-experiment}

Let $\Theta\sim\Pi$ on $(\mathsf{\Theta},\mathcal A)$, and let $X^{(n)}\in(\mathcal X_n,\mathcal B_n)$ have conditional law $P_\theta^{(n)}$ given $\Theta=\theta$. The joint law of the Bayesian experiment is
\[
P_{\Theta,X^{(n)}}(d\theta,dx)=\Pi(d\theta)P_\theta^{(n)}(dx),
\]
where $n$ usually denotes sample size, but may represent an information level or inverse noise scale in sequence and inverse problems.

Whenever finite, the parameter--data mutual information is
\begin{equation*}
%	\label{eq:mutual-information}
	I(\Theta;X^{(n)})
	=
	\KL\!\left(
	P_{\Theta,X^{(n)}}\,\middle\|\,\Pi\otimes P_{X^{(n)}}
	\right).
\end{equation*}
If densities exist, $I(\Theta;X^{(n)})=\E[\log\{p(X^{(n)}\mid\Theta)/p(X^{(n)})\}]$. Equivalently,
\begin{equation*}
%	\label{eq:mi-as-expected-kl}
	I(\Theta;X^{(n)})
	=
	\E_{X^{(n)}}
	\KL\{\Pi(\cdot\mid X^{(n)})\|\Pi\}.
\end{equation*}
Thus mutual information is the expected posterior-to-prior Kullback--Leibler divergence under the prior-predictive law. It is a property of the full prior--sampling experiment, not of one realized posterior. Moreover, $I(\Theta;X^{(n)})\geq0$, with equality exactly when $\Theta$ and $X^{(n)}$ are independent.

\subsection{Regular-parametric calibration}
\label{subsec:calibration}

The normalization used below is calibrated by the information growth of a regular finite-dimensional experiment. The following standard sufficient conditions make this benchmark explicit.

\begin{proposition}[Regular-parametric information growth]
	\label{prop:regular-calibration}
	Let $\Theta$ take values in an open set $\mathsf{\Theta}\subset\mathbb R^d$, and suppose $X_1,\ldots,X_n$ are conditionally i.i.d. given $\Theta=\theta$, with density $p_\theta$. Assume that the prior is supported on a compact set $K\subset\operatorname{int}(\mathsf{\Theta})$ and:
	
	\begin{enumerate}
		\item[(A1)] the model is identifiable, the log likelihood is three times continuously differentiable in a neighborhood of $K$, and the derivatives required for the local quadratic expansion are uniformly integrable over $K$;
		
		\item[(A2)] the Fisher information matrix $J(\theta)$ is continuous on $K$ and satisfies $0<m\leq\lambda_{\min}\{J(\theta)\}\leq\lambda_{\max}\{J(\theta)\}\leq M<\infty$ uniformly in $\theta\in K$;
		
		\item[(A3)] \(K\Subset\mathsf{\Theta}\) has regular boundary, and the restriction of the prior density to \(K\) admits a continuously differentiable, strictly positive extension to a neighborhood of \(K\), with
		\[
		0<\underline\pi\leq\pi(\theta)\leq\overline\pi<\infty,
		\qquad \theta\in K.
		\]
			
		\item[(A4)] the likelihood satisfies a uniform local asymptotic normality expansion on $K$, and the corresponding Bernstein--von Mises approximation holds uniformly in $\theta\in K$.
	\end{enumerate}
	
	Then
	\begin{equation}
		\label{eq:regular-mi-expansion}
		I(\Theta;X^{(n)})
		=
		\frac d2\log n+C_\Pi+o(1),
	\end{equation}
	where
	\[
	C_\Pi
	=
	-\frac d2\log(2\pi e)
	+
	\E_\Pi\!\left[
	\frac12\log\det J(\Theta)-\log\pi(\Theta)
	\right].
	\]
	In particular, $I(\Theta;X^{(n)})=(d/2)\log n+O(1)$.
\end{proposition}

\begin{proof}
	For each fixed $\theta\in K$, the regular Bayesian information expansion gives
	\[
	\KL\!\left(P_\theta^{(n)}\,\middle\|\,P_{X^{(n)}}\right)
	=
	\frac d2\log\frac{n}{2\pi e}
	+
	\frac12\log\det J(\theta)
	-
	\log\pi(\theta)
	+
	o(1),
	\]
	with the remainder uniform over $\theta\in K$ under (A1)--(A4); see \citet{clarke1990information}. Averaging with respect to $\Pi$ and using
	\[
	I(\Theta;X^{(n)})
	=
	\E_\Pi
	\KL\!\left(P_\Theta^{(n)}\,\middle\|\,P_{X^{(n)}}\right)
	\]
	gives \eqref{eq:regular-mi-expansion}. Compactness of $K$, the uniform spectral bounds on $J(\theta)$, and the upper and lower bounds on $\pi(\theta)$ ensure that the terms entering $C_\Pi$ are bounded and hence integrable.
\end{proof}

The leading term in Proposition~\ref{prop:regular-calibration} gives the calibration used throughout: one regular parameter direction contributes asymptotically $(1/2)\log n$ units of expected parameter--data information.

\begin{definition}[Bayesian effective dimension]
	\label{def:effdim}
	For $n\geq3$ and $I(\Theta;X^{(n)})<\infty$, define
	\begin{equation}
		\label{eq:effdim-definition}
		d_{\rm eff}(n)=\frac{2I(\Theta;X^{(n)})}{\log n}.
	\end{equation}
\end{definition}

For a regular fixed-prior $d$-dimensional model satisfying Proposition~\ref{prop:regular-calibration}, $d_{\rm eff}(n)=d+O(1/\log n)$ and hence $d_{\rm eff}(n)\to d$. Away from this setting, $d_{\rm eff}(n)$ need not be integer-valued, bounded by the ambient dimension, or monotone in $n$. It is therefore best interpreted as expected Bayesian information expressed in regular-parametric dimension units.

The prior is part of the experiment. Different priors under the same likelihood can produce different values of $d_{\rm eff}(n)$ because they leave different amounts of uncertainty for the data to resolve. If the prior depends on $n$, we write $d_{\rm eff}(n;\Pi_n)$ when needed; the fixed-prior result above remains a calibration, not an upper bound.

\begin{remark}[Choice of information index]
	\label{rem:information-index}
	The normalization is attached to the index used to describe the sequence of experiments. Replacing $n$ by $cn$, with fixed $c>0$, leaves the leading logarithmic calibration unchanged, whereas a polynomial re-indexing can rescale the limiting effective dimension. 
	\end{remark}

\subsection{Structural properties}
\label{subsec:basic-properties}

The following properties are inherited from mutual information, but they clarify the behavior of the normalized quantity.

\begin{proposition}[Reparameterization invariance]
	\label{prop:invariance}
	Let $\phi:\mathsf{\Theta}\to\mathsf{\Phi}$ be a measurable bijection with measurable inverse and set $\Phi=\phi(\Theta)$. Then $I(\Phi;X^{(n)})=I(\Theta;X^{(n)})$ and
	\(
	d_{\rm eff}(n;\Phi)=d_{\rm eff}(n;\Theta).
	\)
\end{proposition}

\begin{proof}
	Let
	\(
	T:\mathsf{\Theta}\times\mathcal X_n
	\longrightarrow
	\mathsf{\Phi}\times\mathcal X_n,
	\;
	T(\theta,x)=(\phi(\theta),x).
	\)
	Since \(\phi\) is a measurable bijection with measurable inverse, so is \(T\). Writing \(P_{\Theta,X}\) and \(P_{\Phi,X}\) for the joint laws of \((\Theta,X^{(n)})\) and \((\Phi,X^{(n)})\), respectively, we have
	\(
	P_{\Phi,X}=P_{\Theta,X}\circ T^{-1}.
	\)
	Likewise,
	\(
	P_\Phi\otimes P_X
	=
	(P_\Theta\otimes P_X)\circ T^{-1}.
	\) Since relative entropy is invariant under measurable bijections, we have,
	\[
	\begin{aligned}
		I(\Phi;X^{(n)})
		&=
		\KL\!\left(
		P_{\Phi,X}\,\middle\|\,P_\Phi\otimes P_X
		\right)\\
		&=
		\KL\!\left(
		P_{\Theta,X}\circ T^{-1}
		\,\middle\|\,
		(P_\Theta\otimes P_X)\circ T^{-1}
		\right)\\
		&=
		\KL\!\left(
		P_{\Theta,X}\,\middle\|\,P_\Theta\otimes P_X
		\right)
		=
		I(\Theta;X^{(n)}).
	\end{aligned}
	\]
	Multiplication by \(2/\log n\) gives the corresponding invariance of
	\(d_{\mathrm{eff}}(n)\).
\end{proof}

Thus one-to-one changes of coordinates, including rescaling, leave the effective dimension unchanged.

\begin{proposition}[Additivity for product experiments]
	\label{prop:additivity}
	For $k=1,2$, let $(\Theta_k,X_k^{(n)})$ be Bayesian experiments satisfying $(\Theta_1,X_1^{(n)})\perp(\Theta_2,X_2^{(n)})$. Then
	\[
	I\!\left((\Theta_1,\Theta_2);(X_1^{(n)},X_2^{(n)})\right)
	=
	I(\Theta_1;X_1^{(n)})
	+
	I(\Theta_2;X_2^{(n)}).
	\]
	If the two experiments share the same information index $n$, then
	\[
	d_{\rm eff}^{(1\otimes2)}(n)
	=
	d_{\rm eff}^{(1)}(n)+d_{\rm eff}^{(2)}(n).
	\]
\end{proposition}

\begin{proof}
	By assumption,
	\(
	(\Theta_1,X_1^{(n)})
	\perp
	(\Theta_2,X_2^{(n)}).
	\)
	Consequently the joint law factors as
	\(
	P_{\Theta_1,\Theta_2,X_1,X_2}
	=
	P_{\Theta_1,X_1}\otimes P_{\Theta_2,X_2},
	\)
	while the product of the parameter and data marginals is
	\[
	P_{\Theta_1,\Theta_2}\otimes P_{X_1,X_2}
	=
	(P_{\Theta_1}\otimes P_{X_1})
	\otimes
	(P_{\Theta_2}\otimes P_{X_2}).
	\]
	Additivity of Kullback--Leibler divergence over product measures therefore gives
	\[
	\begin{aligned}
		I\!\left(
		(\Theta_1,\Theta_2);
		(X_1^{(n)},X_2^{(n)})
		\right)
		&=
		\KL\!\left(
		P_{\Theta_1,X_1}\otimes P_{\Theta_2,X_2}
		\,\middle\|\,
		(P_{\Theta_1}\otimes P_{X_1})
		\otimes
		(P_{\Theta_2}\otimes P_{X_2})
		\right)\\
		&=
		I(\Theta_1;X_1^{(n)})
		+
		I(\Theta_2;X_2^{(n)}).
	\end{aligned}
	\]
	Since both experiments use the same information index \(n\), division by
	\((1/2)\log n\) yields
	\(
	d_{\mathrm{eff}}^{(1\otimes2)}(n)
	=
	d_{\mathrm{eff}}^{(1)}(n)
	+
	d_{\mathrm{eff}}^{(2)}(n).
	\)
\end{proof}

The common-index qualification is essential: mutual information is intrinsically additive under independent products, whereas exact additivity of the normalized quantity depends on using the same denominator.

For any statistic or randomized summary $W$ constructed at information level $n$, write $d_{\rm eff}(W;n)=2I(\Theta;W)/\log n$.

\begin{proposition}[Data processing]
	\label{prop:data-processing}
	If $\Theta\to X^{(n)}\to Y\to Z$ is a Markov chain, then
	\[
	I(\Theta;Z)\leq I(\Theta;Y)\leq I(\Theta;X^{(n)}),
	\]
	and hence
	\[
	d_{\rm eff}(Z;n)
	\leq
	d_{\rm eff}(Y;n)
	\leq
	d_{\rm eff}(X^{(n)};n).
	\]
\end{proposition}

\begin{proof}
	The assumed Markov chain
	\(
	\Theta\longrightarrow X^{(n)}\longrightarrow Y\longrightarrow Z
	\)
	allows the data-processing inequality to be applied twice. First,
	\(
	I(\Theta;Y)\leq I(\Theta;X^{(n)}),
	\)
	and then
	\(
	I(\Theta;Z)\leq I(\Theta;Y).
	\)
	Thus
	\(
	I(\Theta;Z)
	\leq
	I(\Theta;Y)
	\leq
	I(\Theta;X^{(n)}).
	\)
	Because \(\log n>0\) for \(n\geq3\), multiplication by \(2/\log n\) preserves these inequalities and gives
	\(
	d_{\mathrm{eff}}(Z;n)
	\leq
	d_{\mathrm{eff}}(Y;n)
	\leq
	d_{\mathrm{eff}}(X^{(n)};n).
	\)
\end{proof}

This includes deterministic statistics, randomized summaries, and other forms of coarsening. A sufficient statistic preserves mutual information; information discarded by further processing cannot be recovered by the normalized scale.

Taken together, reparameterization invariance, product additivity under a common information index, and data processing make $d_{\rm eff}(n)$ a structurally useful information scale. These properties arise from mutual information itself and are not specific to the particular normalization.

\subsection{Relation to other complexity measures}
\label{subsec:relation-complexity}

The primitive quantity in \eqref{eq:effdim-definition} is classical expected information gain. The present construction uses its growth relative to the regular-parametric $(1/2)\log n$ benchmark. This distinction is important: parameter--data mutual information, posterior-to-prior information, and dependence-based Bayesian complexity all have established precedents, while the role of $d_{\rm eff}(n)$ here is to place that information on a sample-size-dependent dimension scale.

The resulting quantity targets a different aspect of complexity from posterior contraction, estimator degrees of freedom, singular learning coefficients, and operator-based intrinsic dimensions, each of which describes a different statistical or computational object. The normalization is particularly informative when parameter–data information is distributed unevenly across directions, as the Gaussian, inverse-problem, and hierarchical examples below illustrate.

\section{Gaussian examples and spectral structure}
\label{sec:gaussian-examples}

Gaussian models make the proposed normalization particularly transparent. Mutual information can be calculated exactly, and its dependence on the prior and likelihood is carried by a spectrum of signal-to-noise ratios. The calculations in this section are finite-sample identities. They also give the simplest examples in which \(d_{\mathrm{eff}}(n)\) can differ appreciably from the ambient parameter dimension.

\subsection{A Gaussian channel identity}
\label{subsec:gaussian-channel}

We begin with the standard linear Gaussian channel in a form that will be used repeatedly.

\begin{lemma}[Gaussian mutual information]
	\label{lem:gaussian-channel-mi}
	Let \(\Theta\sim N(0,\Sigma_\Theta)\) and \(\varepsilon\sim N(0,\Sigma_\varepsilon)\) be independent, with \(\Theta\in\mathbb R^p\), \(\varepsilon\in\mathbb R^m\), \(\Sigma_\Theta\succeq0\), and \(\Sigma_\varepsilon\succ0\). For \(Y=A\Theta+\varepsilon\), \(A\in\mathbb R^{m\times p}\),
	\begin{equation}
		\label{eq:mi-linear-gaussian}
		I(\Theta;Y)
		=
		\frac12\log\det\!\left(
		I_m+\Sigma_\varepsilon^{-1/2}A\Sigma_\Theta A^\top\Sigma_\varepsilon^{-1/2}
		\right).
	\end{equation}
	Equivalently,
	\[
	I(\Theta;Y)
	=
	\frac12\log\det\!\left(
	I_p+\Sigma_\Theta^{1/2}A^\top\Sigma_\varepsilon^{-1}A\Sigma_\Theta^{1/2}
	\right).
	\]
	If \(\rho_1,\ldots,\rho_r>0\) are the nonzero eigenvalues of
	\(\Sigma_\varepsilon^{-1/2}A\Sigma_\Theta A^\top\Sigma_\varepsilon^{-1/2}\), then
	\begin{equation}
		\label{eq:gaussian-mi-spectral}
		I(\Theta;Y)=\frac12\sum_{j=1}^r\log(1+\rho_j).
	\end{equation}
\end{lemma}

\begin{proof}
	The marginal law of \(Y\) is Gaussian with covariance \(\Sigma_Y=A\Sigma_\Theta A^\top+\Sigma_\varepsilon\). Hence \(I(\Theta;Y)=h(Y)-h(Y\mid\Theta)=\frac12\log\{\det(\Sigma_Y)/\det(\Sigma_\varepsilon)\}\), which gives \eqref{eq:mi-linear-gaussian}. The parameter-space form follows from Sylvester's determinant identity, and diagonalizing the positive semidefinite matrix in \eqref{eq:mi-linear-gaussian} gives \eqref{eq:gaussian-mi-spectral}.
\end{proof}

The spectral representation is the useful part. Each nonzero eigendirection contributes \((1/2)\log(1+\rho_j)\), with \(\rho_j\) acting as a directional signal-to-noise ratio. Weak directions contribute approximately \(\rho_j/2\), whereas strong directions continue to accumulate information only logarithmically.

\subsection{Gaussian location models}
\label{subsec:location}

Consider first the scalar model \(X_i\mid\theta\stackrel{\mathrm{iid}}{\sim}N(\theta,\sigma^2)\), \(\theta\sim N(0,\tau^2)\). Since \(\bar X\mid\theta\sim N(\theta,\sigma^2/n)\), Lemma~\ref{lem:gaussian-channel-mi} gives
\[
I(\theta;X_{1:n})=\frac12\log\left(1+\frac{n\tau^2}{\sigma^2}\right),
\qquad
d_{\mathrm{eff}}(n)=\frac{\log(1+n\tau^2/\sigma^2)}{\log n}\to1.
\]

The same calculation in \(d\) dimensions, with \(\theta\sim N_d(0,\tau^2I_d)\) and \(X_i\mid\theta\stackrel{\mathrm{iid}}{\sim}N_d(\theta,\sigma^2I_d)\), yields
\[
I(\theta;X_{1:n})=\frac d2\log\left(1+\frac{n\tau^2}{\sigma^2}\right),
\qquad
d_{\mathrm{eff}}(n)
=
d\,\frac{\log(1+n\tau^2/\sigma^2)}{\log n}
\to d.
\]

This elementary example already shows why the prior belongs in the definition. If \(\tau^2=\tau_n^2\), then
\[
d_{\mathrm{eff}}(n;\Pi_n)
=
d\,\frac{\log(1+n\tau_n^2/\sigma^2)}{\log n}.
\]
Thus \(\tau_n^2\asymp n^{-\gamma}\) gives \(d_{\mathrm{eff}}(n;\Pi_n)\to d(1-\gamma)_+\), while \(\tau_n^2\asymp n^\gamma\), \(\gamma>0\), gives \(d_{\mathrm{eff}}(n;\Pi_n)\to d(1+\gamma)\). The ordinary dimension \(d\) is recovered under the fixed-prior calibration. If the prior itself changes at a polynomial rate, the normalized information records that change as well.

\subsection{Gaussian linear regression}
\label{subsec:linear-regression}

Now consider
\begin{equation}
	\label{eq:linreg}
	Y\mid\beta\sim N(X\beta,\sigma^2I_n),
	\qquad
	\beta\sim N(0,\tau^2I_p),
\end{equation}
where \(X\in\mathbb R^{n\times p}\) is deterministic. Write \(\lambda=\tau^2/\sigma^2\), \(r=\operatorname{rank}(X)\), and let \(s_1\ge\cdots\ge s_r>0\) be the nonzero singular values of \(X\).

\begin{proposition}[Spectral form in Gaussian regression]
	\label{prop:linreg-mi}
	Under \eqref{eq:linreg},
	\begin{equation*}
%		\label{eq:linreg-mi}
		I(\beta;Y)
		=
		\frac12\log\det(I_p+\lambda X^\top X)
		=
		\frac12\sum_{j=1}^r\log(1+\lambda s_j^2),
	\end{equation*}
	and consequently
	\begin{equation}
		\label{eq:linreg-deff}
		d_{\mathrm{eff}}(n)
		=
		\frac1{\log n}\sum_{j=1}^r\log(1+\lambda s_j^2).
	\end{equation}
\end{proposition}

\begin{proof}
	Apply Lemma~\ref{lem:gaussian-channel-mi} with \(A=X\), \(\Sigma_\Theta=\tau^2I_p\), and \(\Sigma_\varepsilon=\sigma^2I_n\). The nonzero eigenvalues of \(X^\top X\) are \(s_1^2,\ldots,s_r^2\).
\end{proof}

Formula \eqref{eq:linreg-deff} shows where dimensional reduction enters. A coefficient direction contributes only through the corresponding singular mode of the design. Directions outside the range of \(X\) contribute nothing, and weak singular modes may contribute very little even when \(p\) is large. Since \(d_{\mathrm{eff}}(n)\le r\log(1+\lambda s_1^2)/\log n\), we have \(d_{\mathrm{eff}}(n)=O(r)\) under \(s_1^2=O(n)\) and fixed \(\lambda\). There is no universal inequality \(d_{\mathrm{eff}}(n)\le r\) without further scaling assumptions.

A finer picture appears when the singular values grow at different rates. Suppose \(s_{j,n}^2\asymp n^{\gamma_j}\), with fixed \(r\), \(\gamma_j\ge0\), and fixed \(\lambda>0\). Then \(\log(1+\lambda s_{j,n}^2)/\log n\to\gamma_j\), so
\(
	d_{\mathrm{eff}}(n)\longrightarrow\sum_{j=1}^r\gamma_j.
\)
The normalized-design case corresponds to \(\gamma_j=1\) for each well-identified direction and gives the usual rank \(r\). Modes that accumulate information more slowly contribute fractionally.

\subsection{Connection with ridge degrees of freedom}
\label{subsec:df}

The same Gaussian model gives a direct connection with a familiar frequentist complexity measure. Under \eqref{eq:linreg}, the posterior mean is the ridge estimator
\[
\hat\beta_\alpha=(X^\top X+\alpha I_p)^{-1}X^\top Y,
\qquad
\alpha=\frac{\sigma^2}{\tau^2}=\frac{1}{\lambda}.
\]
The fitted values are \(\hat Y_\alpha=S_\alpha Y\), where
\[
S_\alpha=X(X^\top X+\alpha I_p)^{-1}X^\top,
\]
and the corresponding ridge degrees of freedom are
\[
\operatorname{df}(\alpha)
=
\operatorname{tr}(S_\alpha)
=
\sum_{j=1}^r \frac{s_j^2}{s_j^2+\alpha}.
\]

Connections between Bayesian information and effective parameter counts have appeared previously, including dependence-based constructions in \citet{van2012bayesian}. More generally, Gaussian mutual information is linked to estimation error through the I--MMSE identity of \citet{guo2005mutual}. In the present conjugate regression model, the same structure yields a particularly simple relation with classical ridge degrees of freedom.

Let
\[
\mathcal I(\lambda)
=
I_\lambda(\beta;Y)
=
\frac12\sum_{j=1}^r \log(1+\lambda s_j^2),
\]
where \(\lambda=\tau^2/\sigma^2\) is the prior-to-noise ratio.

\begin{proposition}[Mutual information as integrated ridge degrees of freedom]
	\label{prop:mi-df-identity}
	For every \(\lambda>0\),
	\begin{equation}
		\label{eq:mi-df-derivative}
		2\lambda \mathcal I'(\lambda)
		=
		\operatorname{df}(1/\lambda).
	\end{equation}
	Consequently,
	\begin{equation}
		\label{eq:mi-df-integral}
		2\mathcal I(\lambda)
		=
		\int_0^\lambda
		\frac{\operatorname{df}(1/t)}{t}\,dt.
	\end{equation}
\end{proposition}

\begin{proof}
	Differentiating,
	\[
	2\lambda \mathcal I'(\lambda)
	=
	\sum_{j=1}^r
	\frac{\lambda s_j^2}{1+\lambda s_j^2}
	=
	\sum_{j=1}^r
	\frac{s_j^2}{s_j^2+1/\lambda}
	=
	\operatorname{df}(1/\lambda),
	\]
	which proves \eqref{eq:mi-df-derivative}. Since \(\mathcal I(0)=0\), integration over \(t\in(0,\lambda)\) gives \eqref{eq:mi-df-integral}.
\end{proof}

Equivalently,
\[
\frac{d}{d\log\lambda}\{2\mathcal I(\lambda)\}
=
\operatorname{df}(1/\lambda).
\]
Thus ridge degrees of freedom describe the local change in parameter--data information along the log signal-to-noise scale, while mutual information accumulates this change over lower scales. A strongly identified mode contributes nearly one unit to ridge degrees of freedom, whereas its accumulated mutual information continues to grow logarithmically.

This relation differs from posterior-predictive dependence constructions that recover degrees of freedom directly. Here ordinary parameter--data mutual information is represented as an integral of the ridge degrees-of-freedom path over the prior-to-noise scale.

Finally, using \(u/(1+u)\leq \log(1+u)\leq u\) for \(u\geq0\),
\[
\operatorname{df}(1/\lambda)
\leq
2I(\beta;Y)
\leq
\lambda\,\operatorname{tr}(X^\top X).
\]
These inequalities provide useful finite-sample comparisons, although the differential identity gives the sharper structural relation.

\section{Information growth beyond finite rank}
\label{sec:consequences}

The Gaussian calculations of Section~\ref{sec:gaussian-examples} suggest that, beyond finite rank, information growth is governed more naturally by the spectrum of the prior--likelihood pair than by algebraic dimension. We now develop this viewpoint for infinite-dimensional Gaussian experiments.

In Gaussian regression, \(d_{\rm eff}(n)=(\log n)^{-1}\sum_{j=1}^r\log(1+\lambda s_j^2)\), so directions outside the range of \(X\) contribute nothing, while weak singular modes contribute only a fraction of the regular \(\log n\) scale. The prior enters through the same spectrum: changing \(\lambda=\tau^2/\sigma^2\), or allowing it to depend on \(n\), can alter both finite-sample contributions and their asymptotic order. Effective dimension is therefore attached to a sequence of prior--sampling experiments, not to the parameter space alone.

\subsection{Gaussian sequence models}
\label{subsec:gaussian-sequence}

Consider the infinite-dimensional Gaussian experiment
\[
Y_j=\kappa_j\theta_j+n^{-1/2}\varepsilon_j,\qquad
\varepsilon_j\stackrel{\rm ind}{\sim}N(0,\sigma^2),\qquad
\theta_j\stackrel{\rm ind}{\sim}N(0,\tau_j^2),\qquad j\geq1,
\]
where \(\kappa_j>0\) are the singular values of a forward operator. Define the effective signal-to-noise coefficients \(a_j=\kappa_j^2\tau_j^2/\sigma^2\), and, without loss of generality, arrange them so that \(a_1\geq a_2\geq\cdots>0\). Whenever \(\sum_j a_j<\infty\), Lemma~\ref{lem:gaussian-channel-mi} applied coordinatewise gives \(2I(\theta;Y)=\sum_{j\geq1}\log(1+na_j)\) and
\[
d_{\rm eff}(n)=\frac{1}{\log n}\sum_{j\geq1}\log(1+na_j).
\]
Thus the forward operator and the prior enter only through the combined spectrum \((a_j)\).

A coordinate-free summary of this spectrum is its counting function \(N(t)=\#\{j:a_j\geq t\}\), \(t>0\). Since \(\log(1+na_j)=\int_0^{na_j}(1+t)^{-1}dt\), Tonelli's theorem gives the exact layer-cake representation
\begin{equation}
	\label{eq:counting-representation}
	2I(\theta;Y)=\int_0^\infty\frac{N(t/n)}{1+t}\,dt.
\end{equation}
This identity links information accumulation to the number of spectral directions exceeding each signal-to-noise level. The next result is a Tauberian-type statement for this logarithmic spectral functional.

\begin{theorem}[Effective dimension from spectral counting]
	\label{thm:counting-function}
	Suppose that \(N(t)\sim C\,t^{-\alpha}L(1/t)\) as \(t\downarrow0\), where \(C>0\), \(0<\alpha<1\), and \(L\) is slowly varying at infinity. Then \(\sum_{j\geq1}a_j<\infty\), so the information sum is finite for every \(n\), and
	\[
	2I(\theta;Y)\sim C\pi\csc(\pi\alpha)\,n^\alpha L(n).
	\]
	Consequently,
	\[
	d_{\rm eff}(n)\sim C\pi\csc(\pi\alpha)\,\frac{n^\alpha L(n)}{\log n}.
	\]
\end{theorem}

\begin{proof}
	First, \(\sum_{j\geq1}a_j=\int_0^{a_1}N(u)\,du\). Since \(N(u)\sim Cu^{-\alpha}L(1/u)\) and \(0<\alpha<1\), regular variation implies that this integral is finite near zero. Hence \(\sum_j a_j<\infty\), and \(\sum_j\log(1+na_j)\leq n\sum_j a_j<\infty\).
	
	Set \(R(x)=N(1/x)\). Then \(R\) is regularly varying at infinity with index \(\alpha\), and \(R(n)\sim Cn^\alpha L(n)\). By \eqref{eq:counting-representation},
	\[
	\frac{2I(\theta;Y)}{R(n)}
	=
	\int_0^\infty
	\frac{R(n/t)}{R(n)}
	\frac{dt}{1+t},
	\]
	with the convention that \(R(n/t)=N(t/n)\).
	
	Fix \(\varepsilon\in(0,\min\{\alpha,1-\alpha\})\). Choose \(\delta>0\) sufficiently small that Potter bounds (see, for example, \cite{bingham1987regular}) apply whenever both arguments exceed \(1/\delta\). Split
	\[
	2I(\theta;Y)
	=
	\int_0^{\delta n}\frac{N(t/n)}{1+t}\,dt
	+
	\int_{\delta n}^{na_1}\frac{N(t/n)}{1+t}\,dt,
	\]
	where the integrand vanishes for \(t>na_1\).
	
	On \(0<t\leq\delta n\), both \(n\) and \(n/t\) lie in the asymptotic region for all sufficiently large \(n\). Potter bounds therefore give
	\[
	\frac{R(n/t)}{R(n)}
	\leq
	C_\varepsilon
	\begin{cases}
		t^{-\alpha-\varepsilon}, & 0<t\leq1,\\[1mm]
		t^{-\alpha+\varepsilon}, & 1<t\leq\delta n.
	\end{cases}
	\]
	After multiplication by \((1+t)^{-1}\), the right-hand side is integrable on \((0,\infty)\), because \(\alpha+\varepsilon<1\) and \(\alpha-\varepsilon>0\). Moreover, for every fixed \(t>0\), \(R(n/t)/R(n)\to t^{-\alpha}\). Dominated convergence yields
	\[
	\frac{1}{R(n)}
	\int_0^{\delta n}\frac{N(t/n)}{1+t}\,dt
	\longrightarrow
	\int_0^\infty\frac{t^{-\alpha}}{1+t}\,dt.
	\]
	
	For the remaining part,
	\[
	\int_{\delta n}^{na_1}\frac{N(t/n)}{1+t}\,dt
	\leq
	N(\delta)\log\frac{1+na_1}{1+\delta n}
	=
	O(1).
	\]
	Since \(R(n)\to\infty\), this contribution is \(o\{R(n)\}\). Hence
	\[
	\frac{2I(\theta;Y)}{R(n)}
	\longrightarrow
	\int_0^\infty\frac{t^{-\alpha}}{1+t}\,dt
	=
	\pi\csc(\pi\alpha),
	\]
	where the last identity is the beta-integral formula. Using \(R(n)\sim Cn^\alpha L(n)\) proves the result.
\end{proof}

The theorem shows that the small-eigenvalue abundance of the experiment determines its leading information-growth law. The exponent governing \(N(t)\) becomes the exponent governing mutual information, while calibration by the regular-parametric \(\log n\) scale contributes the additional factor \(1/\log n\). Thus the small-eigenvalue behavior of the prior--likelihood spectrum directly determines the growth of Bayesian information on the regular-parametric dimension scale.

\subsection{Polynomial spectral decay}
\label{subsec:polynomial-spectrum}

Polynomial spectra provide the canonical mildly ill-posed regime.

\begin{corollary}[Polynomial spectral order]
	\label{cor:sequence-growth}
	Suppose that, for some \(q>1\) and constants \(0<c_1\leq c_2<\infty\), \(c_1j^{-q}\leq a_j\leq c_2j^{-q}\), \(j\geq1\). Then \(I(\theta;Y)\asymp n^{1/q}\) and \(d_{\rm eff}(n)\asymp n^{1/q}/\log n\). In particular, if \(\kappa_j^2\asymp j^{-2a}\) and \(\tau_j^2\asymp j^{-1-2s}\), then
	\[
	d_{\rm eff}(n)\asymp\frac{n^{1/(2a+2s+1)}}{\log n}.
	\]
\end{corollary}

\begin{proof}
	The assumed bounds imply
	\[
	\sum_{j\geq1}\log(1+c_1nj^{-q})
	\leq
	2I(\theta;Y)
	\leq
	\sum_{j\geq1}\log(1+c_2nj^{-q}).
	\]
	For the exact power spectrum \(cj^{-q}\), the counting function satisfies \(N_c(t)\sim c^{1/q}t^{-1/q}\). Theorem~\ref{thm:counting-function}, with \(\alpha=1/q\), therefore shows that each bounding sum is of order \(n^{1/q}\). The first assertion follows by comparison. For the mildly ill-posed specification, \(a_j\asymp j^{-(2a+2s+1)}\), so take \(q=2a+2s+1\).
\end{proof}

Exact polynomial decay yields the leading constant.

\begin{corollary}[Sharp asymptotic under polynomial decay]
	\label{cor:sequence-sharp}
	If \(a_j\sim cj^{-q}\), \(q>1\), \(c>0\), then
	\[
	2I(\theta;Y)
	\sim
	\pi c^{1/q}\csc\!\left(\frac{\pi}{q}\right)n^{1/q},
	\]
	and
	\[
	d_{\rm eff}(n)
	\sim
	\pi c^{1/q}\csc\!\left(\frac{\pi}{q}\right)
	\frac{n^{1/q}}{\log n}.
	\]
\end{corollary}

\begin{proof}
	Since \(a_j\sim cj^{-q}\), \(N(t)\sim c^{1/q}t^{-1/q}\) as \(t\downarrow0\). Apply Theorem~\ref{thm:counting-function} with \(\alpha=1/q\), \(C=c^{1/q}\), and \(L\equiv1\).
\end{proof}

A standard mildly ill-posed inverse problem has \(\kappa_j^2\asymp j^{-2a}\) and \(\tau_j^2\asymp j^{-1-2s}\), so the exponent \(1/(2a+2s+1)\) records the joint effect of likelihood and prior geometry. Greater ill-posedness \(a\) or stronger prior smoothing \(s\) slows the growth of information measured in regular-parametric dimension units.

The counting-function formulation also separates this quantity from a hard spectral cutoff. Under polynomial decay there are \(O(n^{1/q})\) modes with \(na_j\gtrsim1\), but their information contributions are unequal: modes near the cutoff contribute \(O(1)\), while the strongest modes contribute on the \(\log n\) scale. Effective dimension aggregates these contributions before calibration.

\subsection{Exponential spectral decay}
\label{subsec:exponential-spectrum}

A qualitatively different regime arises when the effective spectrum decays exponentially. Such behavior is characteristic of severely ill-posed inverse problems. The spectral counting function then grows only logarithmically near zero, corresponding formally to the boundary \(\alpha=0\) of Theorem~\ref{thm:counting-function}. A separate argument is therefore needed.

It is useful to formulate the condition on the logarithmic scale. This both weakens exact exponential equivalence and accommodates polynomial or slowly varying factors multiplying the leading exponential term. The assumptions below are understood for the nonincreasing rearrangement of the effective spectrum; equivalently, the proofs apply to any enumeration satisfying the stated coordinatewise bounds.

\begin{theorem}[Exponential spectral decay]
	\label{thm:exponential-spectrum}
	Suppose that \(-\log a_j\sim cj^\beta\), \(c>0\), \(\beta>0\). Then
	\[
	2I(\theta;Y)
	\sim
	\frac{\beta}{\beta+1}c^{-1/\beta}(\log n)^{1+1/\beta},
	\]
	and consequently
	\[
	d_{\rm eff}(n)
	\sim
	\frac{\beta}{\beta+1}c^{-1/\beta}(\log n)^{1/\beta}.
	\]
\end{theorem}

\begin{proof}
	We first consider the exact spectrum \(b_j=e^{-cj^\beta}\). Write \(L=\log n\) and \(m=(L/c)^{1/\beta}\). Then
	\[
	\sum_{j\geq1}\log(1+nb_j)
	=
	\sum_{j\geq1}\log\{1+\exp(L-cj^\beta)\}.
	\]
	
	Fix \(M>1\). For \(j\leq Mm\), \((L-cj^\beta)/L=1-(j/m)^\beta\). Using \(0\leq\log(1+e^z)-z_+\leq\log2\), \(z\in\mathbb R\), we obtain uniformly for \(j\leq Mm\),
	\[
	\frac{1}{L}\log\{1+\exp(L-cj^\beta)\}
	=
	\left\{1-\left(\frac{j}{m}\right)^\beta\right\}_+
	+
	O(L^{-1}).
	\]
	Hence
	\[
	\frac{1}{mL}
	\sum_{j\leq Mm}\log\{1+\exp(L-cj^\beta)\}
	\longrightarrow
	\int_0^M(1-x^\beta)_+\,dx
	=
	\frac{\beta}{\beta+1}.
	\]
	
	For the tail, \(\log(1+u)\leq u\) gives
	\[
	\sum_{j>Mm}\log\{1+\exp(L-cj^\beta)\}
	\leq
	e^L\sum_{j>Mm}e^{-cj^\beta}.
	\]
	Choose \(M_0\) with \(1<M_0<M\). For sufficiently large \(m\),
	\[
	\sum_{j>Mm}e^{-cj^\beta}
	\leq
	\sum_{j>M_0m}e^{-cj^\beta}
	\leq
	e^{-c\lfloor M_0m\rfloor^\beta}
	+
	\int_{\lfloor M_0m\rfloor}^{\infty}e^{-ct^\beta}\,dt.
	\]
	The standard incomplete-gamma estimate gives, for sufficiently large \(x\),
	\[
	\int_x^\infty e^{-ct^\beta}\,dt
	\leq
	C x^{1-\beta}e^{-cx^\beta}.
	\]
	Since \(c(M_0m)^\beta=M_0^\beta L\), multiplication by \(e^L\) yields
	\[
	e^L\sum_{j>Mm}e^{-cj^\beta}
	\leq
	\exp\{-(M_0^\beta-1)L+o(L)\}
	+
	C m^{1-\beta}
	\exp\{-(M_0^\beta-1)L+o(L)\}
	=
	o(mL),
	\]
	because \(M_0^\beta>1\). Thus, for the exact exponential spectrum,
	\[
	\sum_{j\geq1}\log(1+ne^{-cj^\beta})
	\sim
	\frac{\beta}{\beta+1}c^{-1/\beta}L^{1+1/\beta}.
	\]
	
	We now return to the assumed spectrum. For every \(\varepsilon\in(0,c)\), the relation \(-\log a_j\sim cj^\beta\) implies that, for all sufficiently large \(j\),
	\[
	e^{-(c+\varepsilon)j^\beta}
	\leq
	a_j
	\leq
	e^{-(c-\varepsilon)j^\beta}.
	\]
	The finitely many exceptional coordinates contribute \(O(\log n)=o\!\left((\log n)^{1+1/\beta}\right)\). Comparison with the two exact exponential spectra therefore gives
	\[
	\frac{\beta}{\beta+1}(c+\varepsilon)^{-1/\beta}
	\leq
	\liminf_{n\to\infty}
	\frac{2I(\theta;Y)}{(\log n)^{1+1/\beta}}
	\]
	and
	\[
	\limsup_{n\to\infty}
	\frac{2I(\theta;Y)}{(\log n)^{1+1/\beta}}
	\leq
	\frac{\beta}{\beta+1}(c-\varepsilon)^{-1/\beta}.
	\]
	Letting \(\varepsilon\downarrow0\) proves the first assertion. Division by \(\log n\) gives the result for \(d_{\rm eff}(n)\).
\end{proof}

The polynomial and exponential regimes therefore have qualitatively different information-growth laws: \(a_j\asymp j^{-q}\Longrightarrow d_{\rm eff}(n)\asymp n^{1/q}/\log n\), whereas \(-\log a_j\asymp j^\beta\Longrightarrow d_{\rm eff}(n)\asymp(\log n)^{1/\beta}\). This distinction parallels the usual separation between mildly and severely ill-posed inverse problems, but here it is expressed through prior-to-posterior information growth.

The logarithmic formulation of Theorem~\ref{thm:exponential-spectrum} also covers spectra with polynomial prefactors.

\begin{corollary}[Polynomially modified exponential spectra]
	\label{cor:poly-exponential}
	Suppose that \(a_j=j^{-r}\ell(j)e^{-cj^\beta}\), \(c>0\), \(\beta>0\), where \(r\in\mathbb R\), \(\ell(j)>0\), and \(\log\ell(j)=o(j^\beta)\). Then
	\[
	2I(\theta;Y)
	\sim
	\frac{\beta}{\beta+1}c^{-1/\beta}(\log n)^{1+1/\beta},
	\]
	and
	\[
	d_{\rm eff}(n)
	\sim
	\frac{\beta}{\beta+1}c^{-1/\beta}(\log n)^{1/\beta}.
	\]
\end{corollary}

\begin{proof}
	Since \(-\log a_j=cj^\beta+r\log j-\log\ell(j)=cj^\beta+o(j^\beta)\), the conclusion follows directly from Theorem~\ref{thm:exponential-spectrum}.
\end{proof}

In particular, suppose \(\kappa_j^2\asymp e^{-2aj^\beta}\), \(\tau_j^2\asymp j^{-1-2s}\), \(a>0\), \(s>0\). Then
\[
-\log\left(\frac{\kappa_j^2\tau_j^2}{\sigma^2}\right)
=
2aj^\beta+O(\log j),
\]
and hence
\[
d_{\rm eff}(n)
\sim
\frac{\beta}{\beta+1}(2a)^{-1/\beta}(\log n)^{1/\beta}.
\]
Thus prior smoothness contributes only a lower-order logarithmic correction to the spectral cutoff in this severely ill-posed regime; the leading information-growth scale is determined by the exponential decay of the forward operator. This contrasts with polynomial ill-posedness, where prior smoothness and operator decay jointly determine the leading exponent.

\subsection{What the information scale captures}
\label{subsec:interpretation}

The quantity \(d_{\rm eff}(n)\) summarizes expected prior-to-posterior information in regular-parametric dimension units. In spectral models it automatically discounts weakly identified directions, but it should not be interpreted literally as a count of identifiable coordinates. Nor must it increase with \(n\): for nested observations, \(I(\Theta;X^{(n+1)})=I(\Theta;X^{(n)})+I(\Theta;X_{n+1}\mid X^{(n)})\), so mutual information is nondecreasing, while division by \(\log n\) can reverse finite-sample comparisons. The asymptotic growth of \(d_{\rm eff}(n)\) is therefore more informative than monotonicity at individual information levels.

The sequence-model results also clarify how regularization acts on information growth. Prior smoothing and attenuation by an ill-posed operator modify the effective signal-to-noise spectrum: Theorem~\ref{thm:counting-function} shows that, under regular variation, information growth is governed by how densely this spectrum accumulates near zero, while Theorem~\ref{thm:exponential-spectrum} gives the corresponding regime under exponential decay. Section~\ref{sec:shrinkage} extends this viewpoint to hierarchical models in which the relevant scales are themselves random.

\section{Scale-mixture and shrinkage priors}
\label{sec:shrinkage}

Global--local and scale-mixture priors make the prior geometry itself random. Conditional on the latent scales, the model is often Gaussian and the calculations of Section~\ref{sec:gaussian-examples} apply directly. Marginally, however, the data may also carry information about the latent scales. The distinction between these two contributions is essential: control of the conditional Gaussian term does not by itself imply control of the full marginal mutual information.

\subsection{Scalar scale mixtures}
\label{subsec:gl-priors}

Consider \(Y\mid\theta\sim N(\theta,\sigma^2/n)\), \(\theta\mid\lambda\sim N(0,\lambda^2)\), and \(\lambda\sim\Pi_\lambda\), with \(\lambda>0\). The hierarchy induces the Markov chain \(\lambda\longrightarrow\theta\longrightarrow Y\). By the chain rule, \(I(\theta,\lambda;Y)=I(\lambda;Y)+I(\theta;Y\mid\lambda)\), while the Markov property gives \(I(\lambda;Y\mid\theta)=0\), and therefore \(I(\theta,\lambda;Y)=I(\theta;Y)\). Hence
\begin{equation}
	\label{eq:shrinkage-mi-decomposition}
	I(\theta;Y)=I(\lambda;Y)+I(\theta;Y\mid\lambda).
\end{equation}

Conditional on \(\lambda=\ell\), \(I(\theta;Y\mid\lambda=\ell)=\frac12\log(1+n\ell^2/\sigma^2)\), so
\begin{equation}
	\label{eq:conditional-shrinkage-mi}
	I(\theta;Y)
	=
	I(\lambda;Y)
	+
	\frac12\E_\lambda
	\log\left(1+\frac{n\lambda^2}{\sigma^2}\right).
\end{equation}

Equation~\eqref{eq:conditional-shrinkage-mi} separates two sources of information. The second term is the average Gaussian information available if the scale were observed. The first measures how much the data reveal about the latent scale itself. Any statement about the marginal effective dimension must account for both terms.

\subsection{Finite-moment mixtures}

Let \(c=n/\sigma^2\). If \(\E(\lambda^2)<\infty\), concavity gives \(\E_\lambda\log(1+c\lambda^2)\leq\log\{1+c\,\E(\lambda^2)\}\). Therefore \(I(\theta;Y\mid\lambda)\leq\frac12\log\{1+c\,\E(\lambda^2)\}\). Equivalently, \eqref{eq:conditional-shrinkage-mi} yields
\[
I(\theta;Y)
\leq
I(\lambda;Y)
+
\frac12\log\{1+c\,\E(\lambda^2)\}.
\]
Thus a finite second moment controls only the conditional Gaussian contribution. The latent-scale term \(I(\lambda;Y)\) remains separate.

A centered Student-\(t_\nu\) prior illustrates the distinction. Writing \(\theta\mid\lambda^2\sim N(0,\lambda^2)\) and \(\lambda^2\sim\operatorname{Inv\text{-}Gamma}(\nu/2,\nu s^2/2)\), we have \(\E(\lambda^2)=\nu s^2/(\nu-2)\) for \(\nu>2\), and hence \(I(\theta;Y\mid\lambda)\leq\frac12\log\{1+c\nu s^2/(\nu-2)\}\). When \(\nu\leq2\), this second-moment argument is unavailable, although the conditional Gaussian information may still be finite because it depends only on a logarithmic moment.

\subsection{Heavy-tailed scales}
\label{subsec:heavy-tailed-scales}

The next lemma gives a sufficient condition for finiteness of the conditional Gaussian information even when ordinary scale moments do not exist.

\begin{lemma}[Finite conditional Gaussian information under polynomial tails]
	\label{lem:heavy-tail-log-moment}
	Let \(\lambda\geq0\), and suppose that for some \(\alpha>0\), \(t_0\geq1\), and \(C>0\), \(\Pr(\lambda\geq t)\leq Ct^{-\alpha}\) for \(t\geq t_0\). Then \(\E\log(1+\lambda^2)<\infty\). Moreover, for every \(c>0\),
	\[
	\E\log(1+c\lambda^2)
	\leq
	\log(1+c)
	+
	\log(1+t_0^2)
	+
	\frac{2C}{\alpha}t_0^{-\alpha}.
	\]
	Consequently, \(I(\theta;Y\mid\lambda)=\frac12\E_\lambda\log(1+c\lambda^2)<\infty\), and \(I(\theta;Y\mid\lambda)=O(\log c)\) as \(c\to\infty\).
\end{lemma}

\begin{proof}
	Since \(1+c\lambda^2\leq(1+c)(1+\lambda^2)\), we have \(\log(1+c\lambda^2)\leq\log(1+c)+\log(1+\lambda^2)\). It therefore suffices to control \(\E\log(1+\lambda^2)\). For \(g(t)=\log(1+t^2)\), the tail-integral identity gives
	\[
	\E g(\lambda)
	=
	\int_0^\infty
	\frac{2t}{1+t^2}\Pr(\lambda\geq t)\,dt.
	\]
	The contribution over \([0,t_0]\) is at most \(\log(1+t_0^2)\). Since \(2t/(1+t^2)\leq2/t\) for \(t\geq1\),
	\[
	\int_{t_0}^\infty
	\frac{2t}{1+t^2}\Pr(\lambda\geq t)\,dt
	\leq
	2C\int_{t_0}^\infty t^{-(\alpha+1)}\,dt
	=
	\frac{2C}{\alpha}t_0^{-\alpha}.
	\]
	The stated bounds follow.
\end{proof}

The same tail condition is strong enough to control the full marginal information, despite the fact that \(\lambda\) is continuous.

\begin{proposition}[Finite marginal information under polynomial scale tails]
	\label{prop:scale-mi-finite}
	In the scalar scale-mixture model above, suppose \(\Pr(\lambda\geq t)\leq Ct^{-\alpha}\) for \(t\geq t_0\), with \(\alpha>0\), \(t_0\geq1\), and \(C>0\). Then, for every \(n\geq1\),
	\[
	I(\theta;Y)<\infty.
	\]
	Consequently, both terms in \eqref{eq:shrinkage-mi-decomposition} are finite.
\end{proposition}

\begin{proof}
	Write \(\theta=\lambda Z\), where \(Z\sim N(0,1)\) is independent of \(\lambda\), and \(Y=\theta+\varepsilon=\lambda Z+\varepsilon\), where \(\varepsilon\sim N(0,\sigma^2/n)\) is independent of \((\lambda,Z)\). By Lemma~\ref{lem:heavy-tail-log-moment}, \(\E\log(1+\lambda^2)<\infty\). Since \(1+\lambda^2Z^2\leq(1+\lambda^2)(1+Z^2)\),
	\[
	\E\log(1+\theta^2)
	\leq
	\E\log(1+\lambda^2)+\E\log(1+Z^2)
	<
	\infty.
	\]
	Moreover, \(1+(\theta+\varepsilon)^2\leq2(1+\theta^2)(1+\varepsilon^2)\), and hence
	\[
	\E\log(1+Y^2)
	\leq
	\log2+\E\log(1+\theta^2)+\E\log(1+\varepsilon^2)
	<
	\infty.
	\]
	
	Let \(p_Y\) denote the density of \(Y\) and \(q(y)=\{\pi(1+y^2)\}^{-1}\) the standard Cauchy density. Nonnegativity of \(\KL(p_Y\|q)\) gives
	\[
	h(Y)
	\leq
	-\E\log q(Y)
	=
	\log\pi+\E\log(1+Y^2)
	<
	\infty.
	\]
	On the other hand, \(Y\) is the convolution of the law of \(\theta\) with a nondegenerate Gaussian density, so \(p_Y\) is bounded by \(\sqrt{n/(2\pi\sigma^2)}\). Hence
	\[
	h(Y)\geq-\log\|p_Y\|_\infty>-\infty.
	\]
	Therefore \(h(Y)\) is finite. Since \(Y\mid\theta\sim N(\theta,\sigma^2/n)\),
	\[
	I(\theta;Y)
	=
	h(Y)-h(Y\mid\theta)
	=
	h(Y)-\frac12\log\left(\frac{2\pi e\sigma^2}{n}\right)
	<
	\infty.
	\]
	Finally, \eqref{eq:shrinkage-mi-decomposition} expresses this finite quantity as the sum of two nonnegative terms, so \(I(\lambda;Y)<\infty\) and \(I(\theta;Y\mid\lambda)<\infty\).
\end{proof}

For the horseshoe local scale, \(\lambda\sim C^+(0,1)\), the polynomial-tail condition holds with \(\alpha=1\). Conditional on a fixed global scale \(\tau>0\), replacing \(\lambda\) by \(\tau\lambda\) preserves a polynomial tail of the same order. Lemma~\ref{lem:heavy-tail-log-moment} therefore gives finite conditional Gaussian information for every \(n\), despite \(\E(\lambda^2)=\infty\), while Proposition~\ref{prop:scale-mi-finite} establishes finiteness of the full marginal information \(I(\theta;Y)\). These results concern finiteness; they do not provide sharp asymptotic rates for either component of \eqref{eq:shrinkage-mi-decomposition}.

\subsection{Global--local regression}
\label{subsec:multivariate-shrinkage}

Consider \(Y\mid\beta\sim N(X\beta,\sigma^2I_n)\), \(\beta_j\mid\tau,\lambda_j\sim N(0,\tau^2\lambda_j^2)\), \(j=1,\ldots,p\), and write \(\Sigma_{\lambda,\tau}=\diag(\tau^2\lambda_1^2,\ldots,\tau^2\lambda_p^2)\). Conditional on \((\boldsymbol\lambda,\tau)\), Lemma~\ref{lem:gaussian-channel-mi} gives
\begin{equation*}
%	\label{eq:conditional-mi-regression}
	I(\beta;Y\mid\boldsymbol\lambda,\tau)
	=
	\frac12
	\log\det\left(
	I_n+\sigma^{-2}X\Sigma_{\lambda,\tau}X^\top
	\right).
\end{equation*}
If \(\rho_1,\ldots,\rho_r\) are the nonzero eigenvalues of \(\sigma^{-2}X\Sigma_{\lambda,\tau}X^\top\), then \(I(\beta;Y\mid\boldsymbol\lambda,\tau)=\frac12\sum_{k=1}^r\log(1+\rho_k)\).

The conditional information geometry is therefore determined jointly by the design and the realized shrinkage profile. Even for fixed \(X^\top X\), changing the local scales alters the effective signal-to-noise spectrum seen by the likelihood.

For an orthogonal design, \(X^\top X=\diag(\|X_1\|^2,\ldots,\|X_p\|^2)\), and
\[
I(\beta;Y\mid\boldsymbol\lambda,\tau)
=
\frac12
\sum_{j=1}^p
\log\left(
1+\frac{\tau^2\lambda_j^2\|X_j\|^2}{\sigma^2}
\right).
\]
Small local scales reduce the conditional information contribution of the corresponding coordinates. For a general design there is no coordinatewise decomposition; the scales act through the eigenvalues of \(X\Sigma_{\lambda,\tau}X^\top\).

The full marginal information again contains a separate latent-scale component. If \(L=(\boldsymbol\lambda,\tau)\) collects all scale parameters and the hierarchy satisfies \(L\to\beta\to Y\), then
\begin{equation}
	\label{eq:global-local-mi-decomp}
	I(\beta;Y)
	=
	I(L;Y)
	+
	I(\beta;Y\mid L).
\end{equation}
Thus conditional Gaussian calculations determine only the second term in the marginal information.

This decomposition gives a precise interpretation of global--local shrinkage in information terms. Conditional on the latent scales, shrinkage changes the spectrum governing information accumulation. Marginally, additional information is used to learn the shrinkage configuration itself. Proposition~\ref{prop:scale-mi-finite} shows that polynomially heavy-tailed mixing does not by itself prevent finite marginal information in the scalar Gaussian model. Sharp rates for \(I(L;Y)\) and for the two terms in \eqref{eq:global-local-mi-decomp}, particularly in multivariate models under heavy-tailed global--local priors such as the horseshoe, remain separate theoretical questions.

\section{Open problems and future directions}
\label{sec:open-problems}

The effective dimension studied here suggests several questions that remain open.

A first problem is estimation. Since \(d_{\mathrm{eff}}(n)\) is determined by the joint prior--sampling law, it is not directly observable from a single data set. Practical use would therefore require estimators, bounds, or sensitivity measures that avoid repeated access to the full prior-predictive experiment. Monte Carlo estimates of expected posterior-to-prior Kullback--Leibler divergence, variational representations of mutual information, and local Gaussian approximations provide possible starting points. Understanding their finite-sample accuracy, even in moderately high-dimensional Gaussian models, would be useful.

A second direction concerns experiments beyond the regular Gaussian settings considered here. Singular, weakly identified, misspecified, and non-Gaussian models may exhibit information growth that differs substantially from the regular \((d/2)\log n\) benchmark. In singular models, for example, learning coefficients replace ordinary dimension in marginal-likelihood asymptotics. It would be useful to determine whether analogous quantities govern parameter--data mutual information, and more generally to identify the appropriate information scale when sample size is not the natural index of the experiment.

Finally, hierarchical shrinkage leaves an unresolved information-allocation problem. Section~\ref{sec:shrinkage} separates conditional Gaussian information from information learned about latent scales and establishes marginal finiteness under polynomially heavy-tailed mixing, but the asymptotic sizes of these components remain largely unexplored. Sharp bounds for global--local priors, particularly the horseshoe, could clarify how sparsity and local adaptation are reflected directly in information growth.

These questions would help determine how far the proposed information scale extends beyond the settings where its behavior can currently be characterized explicitly.

\bibliographystyle{apalike}
\bibliography{info-dim-bib}

\end{document}